.
\font\sets=msbm10.
\font\script=eusm10.
\font\stampatello=cmcsc10.
\font\symbols=msam10.
\def\1{{\bf 1}}
\def\defineq{\buildrel{def}\over{=}}
\def\definiz{\buildrel{def}\over{\Longleftrightarrow}}

\def\multisum{\mathop{\sum \cdots \sum}}
\def\C{\hbox{\sets C}}
\def\N{\hbox{\sets N}}
\def\Primes{\hbox{\sets P}}

\def\Z{\hbox{\sets Z}}

\def\cloudF{\hbox{$\buildrel{\subset\!\!\supset}\over{F}$}}
\def\0{{\bf 0}}
\def\cloudzero{$\buildrel{\subset\!\!\supset}\over{\0}$}
\def\Carmichael{{\rm Car}}
\def\Wintner{{\rm Win}}
\def\EssBdd{\hbox{\symbols n}\,}

\def\square{\hbox{\vrule\vbox{\hrule\phantom{s}\hrule}\vrule}}

\par
\centerline{\bf A map of Ramanujan expansions}
\bigskip
\bigskip
\bigskip
\centerline{Giovanni Coppola}
\bigskip
\bigskip
\bigskip
\par
\noindent
{\bf Abstract}. Expanding our talk for NTW2017, the Number Theory Week in Pozna\'{n}, we give a \lq \lq map\rq \rq, for Ramanujan expansions. 
\bigskip
\bigskip
\rightline{\it in honorem of Ramanujan expansions Masters} 
\bigskip
\bigskip
\bigskip
\par
\noindent {\bf 1. Ramanujan Land: Basic definitions, notations and results.}
\enspace 
We call \lq \lq map\rq \rq \thinspace a paper which is neither a survey, nor a paper giving only new results; also, it is not something in which to find complete and state-of-the-art proofs. Simply, it's a tool for people who want to have a very brief and quick look at the main properties regarding that argument: actually, this paper was born because we found ourselves in the situation to learn a great amount of basic facts and shortcuts, regarding the Ramanujan expansions. So, without, of course, aiming at completeness, we embark in the not easy at all task, of supplying a panorama in small scale :  what's a map, if not this? The subject is the very interesting, and still growing, theoretical background of Ramanujan expansions. By the way, the German-style approach (very strong but heavy, for the, say, soft mathematics applied), to these expansions, will be our assumed and quoted knowledge, but not our way to express, here, the main ideas, properties, Lemmas and so on about this argument. Last but not least, so to confirm a map is not a survey, we will also give during the exposition some small results (see our Theorems in $\S3$, for example) of our own, so to speak, handmade production (not to be compared to the huge theory we are trying to embed them in!). We keep our focus and our attention to the new lands we found: $\S5$, $\S6$. In fact, many new results, we found, were inspired by, and came, generalizing $\S6$ results. 

\medskip

\par
We will only give \lq \lq quick\rq \rq \thinspace proofs, always less than two pages, 
otherwise no proof at all (however, we'll try our best to quote good literature in which to find illuminating proofs; by the way, I thank Ram Murty for [M]). 

\bigskip

\par
Prior to other definitions, properties and so on, we give the current definition of Ramanujan expansion.

\medskip

\par
\noindent
Given any $F:\N \rightarrow \C$, we say it has a {\it Ramanujan expansion} (see [CMS] Definition 2) 
$$
F(n)=\sum_{q=1}^{\infty}\widehat{F}(q)c_q(n)
\leqno{(RE)}
$$
\par
\noindent
where the {\it Ramanujan sum} is defined (as Ramanujan himself did in [R], compare [M]), abbreviating hereafter $(a,b)\defineq g.c.d.(a,b)$ the greatest common divisor of $a,b\in \Z$, 
$$
c_q(n)\defineq \sum_{j\le q,(j,q)=1}\cos(2\pi jn/q),
\enspace 
\forall q\in \N, \forall n\in \N, 
\leqno{(RS)}
$$
\par
\noindent
meaning : the series in $(RE)$ {\it converges pointwise} to $F(n)$, $\forall n\in \N$, {\it for certain} $\widehat{F}(q)\in \C$. 
\par
Notice that the definition $(RS)$ may be extended to all $n\in \Z$, thus giving $c_q(0)=\sum_{j\le q,(j,q)=1}1\defineq \varphi(q)$, the {\it Euler function} and using $c_q(-n)=c_q(n)$, from parity of cosine, for negative $-n\in \Z$. 
\smallskip
\par				
The {\bf uniqueness} of these {\it Ramanujan coefficients} \thinspace $\widehat{F}(q)$ \thinspace is {\bf not guaranteed}. A classic example was given by Ramanujan himself, after defining his sums, in [R] for the constant-$0$-arithmetic function, namely \enspace $\0(n)\defineq 0$ $\forall n\in \N$ :  
$$
\sum_{q=1}^{\infty}{1\over q}c_q(n)=\0(n),
\leqno{(C)}
$$
\par
\noindent
thus with non-uniqueness for\footnote{$^1$}{Here the Ram subscript is Ramanujan's \lq \lq copyright\rq \rq \thinspace [R]} \thinspace $\widehat{\0}_{\rm Ram}(q)\defineq 1/q$, as of course also $\widehat{\0}(q)=\0(q)$ gives the same function $\0$. (Notice also: the absolute convergence cannot hold, see the following, in case of $\widehat{\0}_{\rm Ram}$, while for $\widehat{\0}=\0$ it's trivial.) Also, Hardy [H] gave \enspace $\widehat{\0}_{\rm Har}(q)\defineq 1/\varphi(q)$ ! 

\bigskip

The non-uniqueness is the reason why we say that $F$ has \lq \lq a\rq \rq, not \lq \lq the\rq \rq \thinspace Ramanujan expansion of coefficients $\widehat{F}(q)$. The uniqueness of an expansion is not for free, since we need hypotheses on coefficients, as we'll see in $\S3$. 

\bigskip

Expansion (1) is absolutely convergent if, beyond above properties, it has 
$$
\sum_{q=1}^{\infty}\left|\widehat{F}(q)c_q(n)\right|<\infty,
\enspace 
\forall n\in \N. 
$$

\bigskip

\par
To prove that $(C)$ does not converge absolutely, we need a classic Lemma [M], with a \lq \lq very quick\rq \rq, less than half-page, proof. 
\smallskip

\par
\noindent
{\bf Lemma 1.}
{\it Let $q\in \N$. Define as usual $e_q(m)\defineq e^{2\pi im/q}$ $\forall m\in \Z$ the additive characters. Indicate with $\Z_q^*$ the reduced residue classes modulo $q$ and with $\mu$ the M\"{o}bius function [T]. Then } 
$$
c_q(n)=\sum_{j\in \Z_q^*}e_q(jn)=\sum_{{{d|q}\atop {d|n}}}d\mu(q/d)=\varphi(q){{\mu(q/(q,n))}\over {\varphi(q/(q,n))}}. 
$$

\smallskip
\par
\noindent {\bf Proof.} As $q=1$ gives $1$ everywhere and $q=2$ everywhere $(-1)^n$, we assume $q>2$ in the following. The first equality comes from the fact that $j\in \Z_q^* \Leftrightarrow -j\in \Z_q^*$, with $j\not \equiv-j\bmod q$ (from $q>2$), the {\it Euler identity}: 
$$
e_q(jn)=\cos(2\pi jn/q)+i\sin(2\pi jn/q)
$$
\par
\noindent
and the fact that sine function is odd (so the imaginary part vanishes). 
\par
We'll write $\1_{\wp}\defineq 1$ iff (if and only if) $\wp$ is true, $\1_{\wp}\defineq 0$ otherwise, hereafter. 
\par
\noindent
The {\it orthogonality of additive characters}, then g.c.d. rearranging give 
$$
q\1_{q\mid n}=\sum_{r\le q}e_q(rn)
=\sum_{d|q}\sum_{{{r\le q}\atop {(r,q)=q/d}}}e_q(rn)
=\sum_{d|q}\sum_{{{j\le d}\atop {(j,d)=1}}}e_d(jn)
=\sum_{d|q}c_d(n)
\leqno{(1)} 
$$
\par
\noindent
whence the V.I.P. ($=$very important property, like we'll abbreviate hereafter) 
$$
\1_{q|n}={1\over q}\sum_{d|q}c_d(n),
\quad \forall q\in \N, \enspace \forall n\in \Z.
\leqno{(2)} 
$$
\par
\noindent
Then, by M\"{o}bius inversion [T], we get the second equation above: 
$$
q\1_{q|n}=\sum_{d|q}c_d(n)
\enspace \Rightarrow \enspace
c_q(n)=\sum_{d|q}d\1_{d|n}\mu(q/d)=\sum_{{{d|q}\atop {d|n}}}d\mu(q/d). 
$$
\par				
\noindent
This also proves $c_q(n)$ is a multiplicative function of $q$, so we may 
calculate it on the prime-powers $p^K$, as Davenport [Da] does, 
getting the last equation. \hfill $\square$ 

\medskip

This may be applied at once, to prove the next result. 
\smallskip
\par
\noindent
{\bf Proposition 1.}
{\it We have:} 
$$
\sum_{q=1}^{\infty}{1\over q}|c_q(n)|=+\infty,
\quad
\forall n\in \N. 
$$
\par
\noindent {\bf Proof}. In fact, $(q,n)=1$ $\Rightarrow$ $|c_q(n)|=\mu^2(q)$ from Lemma 1, so, fixed $n\in \N$, 
$$
\sum_{q=1}^{\infty}{1\over q}|c_q(n)|\ge \sum_{{{q=1}\atop {(q,n)=1}}}^{\infty}{{\mu^2(q)}\over q}
\ge \sum_{p}{1\over p}-\sum_{p\mid n}{1\over p} = +\infty, 
$$
\par
\noindent
proving the assertion.\hfill $\square$ 

\bigskip

We will not give a list of Ramanujan expansions (even a short one would too long for this paper). See the Ramanujan [R] and Hardy [H] papers, first; then, the book [ScSp] and the two surveys [Lu], [M]. This, in turn, is not an exhaustive list! 

\medskip

\par
\centerline{\bf Notations} 
\par
First, we say $f:\N \rightarrow \C$ satisfies the {\it Ramanujan Conjecture} or, equivalently, it's {\it essentially bounded}, written $f\EssBdd 1$, if $f(n)\ll_{\varepsilon} n^{\varepsilon}$, as $n\to \infty$. Hereafter $\ll$ is the {\it Vinogradov notation}, where $f(n)\ll g(n)$ is equivalent to {\it Landau notation}, $f(n)=O(g(n))$, both meaning that for a certain $n_0\in \N$, we have $|f(n)|\le Cg(n)$, for all $n>n_0$. The constant $C>0$ (named the \lq \lq implicit constant\rq \rq) may depend on other variables, in which case they are displayed as subscripts. Notice, this \lq \lq modified Vinogradov notation\rq \rq, $\EssBdd$, going back at least to Kolesnik in the 50s, doesn't display the $\varepsilon-$dependence explicitly. We write \lq \lq $\ast$\rq \rq, for the {\it Dirichlet product} [T]. The Eratosthenes transform [W], for any $F:\N \rightarrow \C$, is $F'\defineq F\ast \mu$. All of this is standard in Analytic Number Theory: [Da], [T]. 

\medskip 

\par
In the next section we give a kind of candidates, for being Ramanujan coefficients: Wintner's and Carmichael's coefficients. Their main problem is their  existence itself ! However, even if both exist and are equal, they might miss to be Ramanujan coefficients: compare \S7. 
\par
In $\S3$ we describe the existence and also the uniqueness problems, of our Ramanujan coefficients, for general arithmetic functions. 
\par
In $\S4$ we mainly give the two most general theorems for Ramanujan expansions: Hildebrand's [Hi] and Lucht's [Lu]. We shall avoid all the results regarding specifically the additive and the multiplicative arithmetic functions. 
\par
In $\S5$ we turn to a new argument, namely {\bf finite} Ramanujan expansions, abbrev. f.R.e. Here, we explain why do we care about purity,  a property introduced with Ram Murty, in [CM2]; in fact, Hildebrand's result (see $\S4$) gives a kind of f.R.e. for {\bf all} arithmetic functions $F(n)$, with coefficients depending, in general, on $n$: this makes every result and all the calculations very cumbersome, \lq \lq not natural\rq \rq (I am quoting Professor Lutz Lucht, private email conversation) ! 
\par
Our main focus is on $\S5$ and the next $\S6$, in which we introduce [CM2] the {\it shift-Ramanujan expansion}: it was implicit in all of the literature on {\it shifted convolution sums}, expecially for their heuristics. In fact, we are expanding these sums with respect to their shift. We apply here, in finding {\it shift-Ramanujan coefficients}, our results in $\S5$. 
\par
Finally, $\S7$ gives odds and ends, applying shift-Ramanujan expansions to short intervals; also, we show that Carmichael and Wintner coefficients may exist and be the same, without being Ramanujan Coefficients : compare Remark 3 and subsequent V.I.P. 

\bigskip

\par
\noindent {\bf 2. Carmichael, Wintner and Delange Lands: Seeking Ramanujan coefficients.}
\enspace 
We define at once {\it Carmichael's coefficients} for $F:\N \rightarrow \C$ with the property that all the following limits exist in $\C$ : 
$$
\Carmichael_q(F)\defineq {1\over {\varphi(q)}}\lim_{x\to\infty}{1\over x}\sum_{n\le x}F(n)c_q(n), 
\enspace
\forall q\in \N. 
\leqno{(3)}
$$				
\par
\noindent
We define {\it Wintner's coefficients} for $F:\N \rightarrow \C$ with the property that all the following series converge : 
$$
\Wintner_q(F)\defineq \sum_{{{d=1}\atop {d\equiv 0\bmod q}}}^{\infty}{{F'(d)}\over d}, 
\enspace
\forall q\in \N. 
\leqno{(4)}
$$
\par
\noindent
Actually, the {\it Wintner assumption} 
$$
\sum_{d=1}^{\infty}{{|F'(d)|}\over d}<\infty,
\leqno{(WA)} 
$$
\par
\noindent
when satisfied by $F$, ensures at least the existence of all these coefficients in $\C$ (trivially, by positivity). Also, it gives their coincidence with Carmichael's coefficients (which, then, exist), as discovered by Aurel Wintner, [W] : now we call it \lq \lq Wintner's Criterion\rq \rq. 
\smallskip 
\par
\noindent
{\bf Theorem 1.}
\enspace 
{\it Let $F:\N \rightarrow \C$ satisfy Wintner's assumption, namely $(WA)$ above. Then both Carmichael's and Wintner's coefficients exist and they agree.} 
\par
\noindent {\bf Proof}. We prove that Carmichael's coefficients equal Wintner's: for a fixed $q\in \N$ we have 
$$
{1\over {\varphi(q)}}\lim_{x\to \infty}{1\over x}\sum_{n\le x}F(n)c_q(n)=\sum_{d\equiv 0\bmod q}{{F'(d)}\over d}. 
\leqno{(5)}
$$
\par
\noindent
Fix a large $K\in \N$, then plug in LHS (left hand side) the decomposition: 
$$
F(n)=\sum_{d|n,d\le K}F'(d)+\sum_{d|n,d>K}F'(d), 
$$
\par
\noindent
rendering the average in LHS as: 
$$
{1\over x}\sum_{n\le x}F(n)c_q(n)=\sum_{d\le K}F'(d){1\over x}\sum_{m\le x/d}c_q(dm)+\sum_{d>K}F'(d){1\over x}\sum_{m\le x/d}c_q(dm). 
$$
\par
\noindent
We apply two different treatments, depending on $d\le K$ or $d>K$. 
\par
For low divisors $d$, abbreviating \enspace $\Vert \alpha \Vert \defineq \min_{n\in \Z}|\alpha-n|$, we obtain 
$$
\sum_{d\le K}F'(d){1\over x}\sum_{m\le x/d}c_q(dm)=\sum_{d\le K}F'(d)\sum_{j\le q,(j,q)=1}{1\over x}\sum_{m\le x/d}e_q(jdm)
$$
$$
=\sum_{d\le K}F'(d)\sum_{j\le q,(j,q)=1}
 \left({1\over d}\cdot \1_{d\equiv 0\bmod q}
  +O\left({1\over x}\left(1+{{\1_{d\not \equiv 0\bmod q}}\over {\left\Vert jd/q\right\Vert}}\right)\right)\right)
$$
$$
=\varphi(q)\sum_{{{d\le K}\atop {d\equiv 0\bmod q}}}{{F'(d)}\over d}+O(1/x), 
$$
\par
\noindent
from exponential sums cancellations (compare Chapter 25 of [Da]). 
\par
For high divisors $d$, we get  
$$
\sum_{d>K}F'(d){1\over x}\sum_{m\le x/d}c_q(dm)\ll \varphi(q)\sum_{d>K}{{|F'(d)|}\over d}, 
$$
\par
\noindent
uniformly in $x>0$, from the trivial bound $|c_q(n)|\le \varphi(q)$, $\forall n\in \Z$. 
\par				
\noindent
In all, 
$$
{1\over x}\sum_{n\le x}F(n)c_q(n)=\varphi(q)\sum_{{{d\le K}\atop {d\equiv 0\bmod q}}}{{F'(d)}\over d}+O(1/x)+O\Big(\varphi(q)\sum_{d>K}{{|F'(d)|}\over d}\Big), 
$$
\par
\noindent
entailing 
$$
{1\over {\varphi(q)}}\lim_{x\to \infty}{1\over x}\sum_{n\le x}F(n)c_q(n)=\sum_{{{d\le K}\atop {d\equiv 0\bmod q}}}{{F'(d)}\over d}+O\Big(\sum_{d>K}{{|F'(d)|}\over d}\Big), 
$$
\par
\noindent
actually, giving the required equation: from $(WA)$ errors in $O$ are infinitesimal with $K$, an arbitrarily large natural number (also, our LHS doesn't depend on it!). Finally, we get the convergence (to Wintner's coefficients) in RHS (right hand side), of these coefficients, with \thinspace $d\le K$, as $K\to \infty$. \hfill $\square$ 

\medskip

\par
In the literature (compare [ScSp] esp.), the name \lq \lq Wintner's Criterion\rq \rq \thinspace usually refers to the following result (we call Wintner-Delange Formula), optimized in 1976 by Hubert Delange [De]. Theorem 2 explicitly gives an instance of Ramanujan coefficients. 

\medskip

\par
We prove it in this final shape like in [C], taking material from [ScSp] and [De]. See [Lu] for a $4-$lines proof (based on Lucht's Theorem, that we quote, in next section). 
\par 
If \thinspace $\omega(d)\defineq \#\{ p\in \Primes : p|d \}$ \thinspace is the number of prime factors of $d$, then \thinspace $2^{\omega(d)}=\sum_{\ell|d}\mu^2(\ell)$ \thinspace is the number of its square-free [T] divisors. 
\par
We have next {\it Wintner-Delange Formula} for Ramanujan coefficients. 
\smallskip
\par
\noindent
{\bf Theorem 2.}
\enspace 
{\it Let $F:\N \rightarrow \C$ satisfy Delange Hypothesis }
$$
\sum_{d=1}^{\infty}{{2^{\omega(d)}}\over d}\, \left|F'(d)\right|<\infty. 
\leqno{(DH)}
$$
\par
\noindent
{\it Then Wintner's coefficients exist and the Ramanujan expansion}
$$
\sum_{q=1}^{\infty}\Wintner_q(F)c_q(n)
$$ 
\par
\noindent
{\it converges pointwise to $F(n)$, $\forall n\in \N$. Also, Carmichael's coefficients exist and agree with Wintner's\footnote{$^2$}{Of course, $(DH)$ implies $(WA)$, so this part follows from previous result.}.} 
\medskip
\par
\noindent {\bf Proof}. We may confine to the proof of the following formula : 
$$
\sum_{d=1}^{\infty}\sum_{\ell|d}{{F'(d)}\over {d}}c_{\ell}(n)
=\sum_{\ell=1}^{\infty}\sum_{d\equiv 0\bmod \ell}{{F'(d)}\over {d}}c_{\ell}(n), 
\enspace 
\forall n\in \N,
\leqno{(\ast)} 
$$
\par
\noindent
namely, exchanging $\ell,d$ sums in the double series. From this, in fact, $(2)$ gives LHS 
$$
\sum_{d=1}^{\infty}{{F'(d)}\over {d}}\sum_{\ell|d}c_{\ell}(n)
=\sum_{d|n}F'(d)=F(n), 
$$
\par
\noindent
with, in RHS, Wintner's coefficients 
$$
\sum_{d\equiv 0\bmod \ell}{{F'(d)}\over {d}}, 
\enspace \forall \ell\in \N, 
$$
\par
\noindent
thus proving pointwise convergence of Ramanujan expansion in 
$$
F(n)=\sum_{\ell=1}^{\infty}\left(\sum_{d\equiv 0\bmod \ell}{{F'(d)}\over {d}}\right)c_{\ell}(n),
\forall n\in \N. 
$$
\par				
\noindent
In order to prove $(\ast)$, we prove the absolute convergence of the double series involved. 
This, in turn, comes from the fact that LHS of $(\ast)$ with the moduli is 
$$
\sum_{d=1}^{\infty}{{\left|F'(d)\right|}\over {d}}\sum_{\ell|d}\left|c_{\ell}(n)\right|
\le n\sum_{d=1}^{\infty}{{\left|F'(d)\right|}\over {d}}2^{\omega(d)}<\infty, 
\enspace
\forall n\in \N, 
$$
\par
\noindent
that comes from hypothesis $(DH)$, after applying the optimal bound of Hubert Delange: 
$$
\sum_{\ell|d}\left|c_{\ell}(n)\right|\le n\cdot 2^{\omega(d)}, 
$$
\par
\noindent
see [De] (also, for comments about optimality). \hfill $\square$ 

\bigskip

\par
\noindent {\bf 3. Ramanujan clouds and purity: Existence and uniqueness for Ramanujan coefficients.}
\enspace 
From Ramanujan Land (i.e., from his basic definitions) we may see, from the top of its mountain views (i.e., from [R] \& Hardy's [H] expansions of $\0$), the Ramanujan Clouds (i.e., the sets of Ramanujan coefficients for all fixed $F$) ! 
\par
It seems that above Wintner-Delange Formula gives a way to obtain those {\it unique} Ramanujan coefficients, but this is not true. In fact, we may simply add to the Wintner's coefficients (or Carmichael's, which are the same in Delange Hypothesis) the ones of constant-$0$ function; for which we have $\widehat{\0}_{\rm Ram}(q)=1/q$, as we saw, as possible option, like also $\widehat{\0}_{\rm Har}(q)=1/\varphi(q)$. Then, given $F$ with Carmichael's coefficients $\Carmichael_q(F)$, applying Wintner-Delange Formula we get (the same is true with Wintner's coefficients of course) 
$$
\widehat{F}(q)=\Carmichael_q(F)+\alpha/q+\beta/\varphi(q),
\enspace \forall q\in \N, 
$$
\par
\noindent
for all complex numbers $\alpha,\beta$ we like. We see how uniqueness fails! 
\par
The point is that, once we define \lq \lq Ramanujan coefficients\rq \rq, as any sequence such that the series $(RE)$ converges pointwise to $F(n)$, all the \lq \lq Ramanujan coefficients\rq \rq, for the null function, may be added : uniqueness, for the sequence of Ramanujan coefficients, under these hypotheses, is out of sight. However, we have to choose: if we impose too many restrictions we may get no \lq \lq Ramanujan coefficients\rq \rq, at all ! Even worse, this we may not know a priori : we don't know, at the moment, all the conditions to impose for the uniqueness ! 

The task to describe EXPLICITLY ALL of such sequences of coefficients $\hat{F}(q)$, for any given $F:\N \rightarrow \C$, is monumental and amounts exactly to the complete description of the cloud of $\0$. We suspect it has infinite dimension (as a complex vector space). 

Our choice here is to accept non-uniqueness (as Ramanujan did), but we wish to have at least the existence of Ramanujan coefficients. This is guaranteed  by Hildebrand's Theorem : see next section. 

In his survey [Lu], L. Lucht proves in an elementary fashion $(C)$ above and the divisor function $d(a)-$expansion [R] by his  Theorem, see next section. He also discusses non-uniqueness of Ramanujan coefficients. 

\medskip

\par
\noindent {\bf Remark 1.}
\enspace 
{\it In other words, (first of subtleties) the application $\widehat{\enspace}$, sending $F$ to $\widehat{F}$, is not well-defined. On the other hand, the {\it Transforms} $\Carmichael : F \rightarrow \Carmichael(F)$ and resp. $\Wintner : F \rightarrow \Wintner(F)$ giving the sequences of Carmichael, resp., Wintner coefficients, do not exist for all $F$; but, when they exist they are well-defined! We abbreviate $\Carmichael(F)$, resp., $\Wintner(F)$ the sequence $\{ \Carmichael_q(F)\}_{q\in \N}$, resp., $\{ \Wintner_q(F)\}_{q\in \N}$ that are arithmetic functions depending on $F$. } 

\medskip

\par
\noindent
Defining the {\it Ramanujan Cloud} of a given $F:\N \rightarrow \C$, notation \cloudF, is inspired by $F=\0$: 
$$
\hbox{\cloudF}\defineq \left\{ \widehat{F}:\N \rightarrow \C \;\Big|\; \forall n\in \N,\;\sum_{q=1}^{\infty}\widehat{F}(q)c_q(n)\;\hbox{\rm converges pointwise to}\; F(n)\right\} 
$$
\par
\noindent 
and this is always non-empty by Hildebrand's Theorem (see next $\S4$). 
\par 
\noindent
From what we saw above, the cloud of $\0$ contains\enspace ${\widehat{\0}}_{\rm plane}\defineq \{ \alpha\,\widehat{\0}_{\rm Ram} +\beta\,\widehat{\0}_{\rm Har}\,|\,\alpha,\beta\in\C \}$. 
\par				
Given any $F:\N \rightarrow \C$, a so to speak non-trivial property is that, for any fixed sequence of Ramanujan coefficients, say, \thinspace ${\widehat{F}}_0 \defineq \{ {\widehat{F}}_0(q)\}_{q\in\N}$ , the cloud of $F$ contains an important two-dimensional subset 
$$
{\widehat{F}}_0+{\widehat{\0}}_{\rm plane} \subseteq \cloudF 
$$
\par
\noindent
as we saw before. (From Wintner-Delange, may choose ${\widehat{F}}_0=\Carmichael(F)=\Wintner(F)$.) In particular, \enspace ${\widehat{\0}}_{\rm plane}\subseteq$\cloudzero. 

\medskip

Going to open questions (this is 1st): are these inclusions strict ? Better asked, what is the dimension of \cloudzero, as a (complex) vector space ? See that, once we know \cloudzero, all clouds are, trivially, \cloudF$={\widehat{F}}_0+$\cloudzero ! So, how is made \cloudzero ? 

In Ramanujan clouds: {\it given} any {\it two \lq \lq drops\rq \rq}, $\widehat{F}_1,\widehat{F}_2$, {\it in the same cloud}, $\cloudF$ (of a fixed $F:\N \rightarrow \C$), {\it all the line through them}, $\{ \lambda \widehat{F}_1 + (1-\lambda)\widehat{F}_2 : \lambda\in \C\}$, {\it is contained in the same cloud}. Beautiful banner: any two {\it drops see each other in clouds} ! 

\par
For a different approach, based not on hypotheses on the $F$, but on the Ramanujan expansion, we may ask : is there a sufficient condition, ensuring UNIQUENESS of Ramanujan coefficients ? (We understand, from what seen above, that uniqueness, in the hypothesis of pointwise convergence, can't be required!) 
\par
In order to answer, we need first the following definition: 
$$
\sum_{q=1}^{\infty}\widehat{F}(q)c_q(n)
\enspace \hbox{\rm is}\enspace\hbox{\bf pure} \enspace
\enspace \definiz \enspace
\enspace \hbox{\rm all} \enspace \widehat{F}(q) \enspace \hbox{\rm and their supports don't depend on } n. 
$$
\medskip
\par
\noindent
For example, Hildebrand finite R.e.s are not necessarily pure, but $\0$ has a full plane (see above) of pure R.e.s. 
\medskip
\par
From it, we build the definition we need : (hereafter, uniformly convergent$=$converges uniformly in $\N$) 
$$
\sum_{q=1}^{\infty}\widehat{F}(q)c_q(n)
\enspace \hbox{\rm is}\enspace\hbox{\bf completely uniform} \enspace
\enspace \definiz \enspace 
\enspace \hbox{it's pure \& uniformly convergent}
$$
\par
\noindent
that is a kind of strengthening the uniform convergence concept. 
\par
A first UNIQUENESS FORMULA (namely, a formula giving Ramanujan coefficients with uniqueness, under a specific hypothesis) is our Lemma A.4 in [CM2], in which we prove uniqueness for completely uniform Ramanujan expansions: we quote it as stated there. 
\smallskip
\par
\noindent
{\bf Theorem 3.}
\enspace 
{\it Let $F:\N\rightarrow \C$ have an uniformly convergent Ramanujan expansion, i.e. }
$$
F(h)=\sum_{q=1}^{\infty}\widehat{F}(q)c_q(h), 
\enspace 
\forall h\in \N, 
$$
\par
\noindent
{\it with some coefficients $\widehat{F}(q)\in \C$ independent of $h$ (even in their support). Then, these are }
$$
\widehat{F}(\ell)={1\over {\varphi(\ell)}}\lim_{x\to \infty}{1\over x}\sum_{h\le x}F(h)c_{\ell}(h). 
\leqno{(CF)}
$$

\medskip

We call this \lq \lq Carmichael's formula\rq \rq, for Ramanujan coefficients [Ca], [M]. 

\medskip

\par
\noindent {\bf Proof.} Fix $\ell\in \N$ and, by uniform convergence, we have \enspace $\forall \varepsilon>0$ $\exists Q=Q(\varepsilon,\ell)$, with $Q>\ell$ and 
(see soon before Remark 2 for the $d(\ell)$ definition) 
$$
\Big|\sum_{q>Q}\widehat{F}(q)c_q(h)\Big|<{{\varepsilon}\over {d(\ell)}}, 
$$
\par				
\noindent
entailing 
$$
{1\over x}\sum_{h\le x}F(h)c_{\ell}(h)=\sum_{q\le Q}\widehat{F}(q){1\over x}\sum_{h\le x}c_{\ell}(h)c_q(h)+{1\over x}\sum_{h\le x}c_{\ell}(h)\sum_{q>Q}\widehat{F}(q)c_q(h) 
$$
\par
\noindent
(notice purity allows sums exchange) implies (\lq \lq ${\displaystyle \lim_x}$\rq \rq, here, abbreviating \lq \lq ${\displaystyle \lim_{x\to \infty}}$\rq \rq) 
$$
\left|\lim_x {1\over x}\sum_{h\le x}F(h)c_{\ell}(h)
 -\sum_{q\le Q}\widehat{F}(q)\lim_x {1\over x}\sum_{h\le x}c_{\ell}(h)c_q(h)\right|
  \le {{\varepsilon}\over {d(\ell)}} \lim_x {1\over x}\sum_{h\le x}(\ell,h), 
$$
\par
\noindent
from \enspace $|c_{\ell}(h)|\le (\ell,h)$, see Lemma A.1 [CM2], whence the orthogonality relations [M] 
$$
\lim_x {1\over x}\sum_{h\le x}c_{\ell}(h)c_q(h)=\sum_{j\in \Z_{\ell}^*}\sum_{r\in \Z_q^*}\lim_x {1\over x}\sum_{h\le x}e^{2\pi i (j/\ell-r/q)h} 
=\1_{q=\ell}\varphi(\ell)
$$
\par
\noindent
and the formula 
$$
{1\over x}\sum_{h\le x}(\ell,h)=\sum_{t|\ell}{t\over x}\sum_{{{h'\le {x\over t}}\atop {(h',{{\ell}\over t})=1}}}1
=\sum_{t|\ell}{t\over x}\sum_{d\left|{{\ell}\over t}\right.}\mu(d)\left[{x\over {dt}}\right]
=\sum_{t|\ell}\sum_{d\left|{{\ell}\over t}\right.}{{\mu(d)}\over d}+O\left({1\over x}\sum_{t|\ell}td(\ell/t)\right) 
$$
$$
=\sum_{t|\ell}{{\varphi(\ell/t)}\over {\ell/t}}+o(1) 
=\sum_{d|\ell}{{\varphi(d)}\over d}+o(1) 
$$
\par
\noindent
(flipping to the complementary divisor $d=\ell/t$), as $x\to \infty$, together give 
$$
\left|{1\over {\varphi(\ell)}}\lim_{x\to \infty}{1\over x}\sum_{h\le x}F(h)c_{\ell}(h)-\widehat{F}(\ell)\right|
 \le {{\varepsilon}\over {\varphi(\ell)d(\ell)}} \sum_{d|\ell}{{\varphi(d)}\over d}
 \le {{\varepsilon}\over {\varphi(\ell)}}
 \le \varepsilon, 
$$
\par
\noindent
which, as \enspace $\varepsilon>0$ \enspace is arbitrary, shows $(CF)$.\hfill $\square$ 

\bigskip

\par
Our second uniqueness formula is the following new Theorem, the {\bf Wintner-Delange \lq \lq uniqueness\rq \rq \thinspace formula}. See that, instead, Wintner-Delange Formula doesn't give the uniqueness for the Ramanujan coefficients (see the above). 
\smallskip
\par
\noindent
{\bf Theorem 4.}
\enspace 
{\it Let $F:\N \rightarrow \C$ have a pure R.e. }
$$
F(a)=\sum_{q=1}^{\infty}\widehat{F}(q)c_q(a), 
\quad
\forall a\in \N. 
$$
\par
\noindent 
{\it Then}
$$
F'(d)=d\sum_{K=1}^{\infty}\mu(K)\widehat{F}(dK), 
\quad
\forall d\in \N. 
\leqno{(6)}
$$
\par
\noindent 
{\it Furthermore, if the Ramanujan coefficients also satisfy the, say, \lq \lq Dual Delange\rq \rq \thinspace condition }
$$
\sum_{q=1}^{\infty}2^{\omega(q)}\left|\widehat{F}(q)\right|<\infty, 
\leqno{(7)}
$$
\par
\noindent 
{\it then they are unique, as }
$$
\widehat{F}(q)=\Wintner_q(F)=\Carmichael_q(F), 
\quad
\forall q\in \N. 
\leqno{(8)}
$$ 
\smallskip
\par				
\noindent
Of course, both Wintner and Carmichael coefficients in $(8)$ are defined! 
\smallskip
\par
\noindent {\bf Proof.} We start proving $(6)$, inspired by Lucht Theorem, from 
$$
F(n)=\sum_{d|n}d\sum_{K=1}^{\infty}\mu(K)\widehat{F}(dK), 
\leqno{(9)}
$$
\par
\noindent
thanks to Lemma 1, third eq.; then the purity, assumed, implies $(6)$, because: 
RHS of $(6)$ inside RHS of $(9)$ depends only on $d$ (not on $n$ !), so, 
by M\"{o}bius inversion [T], it's $F'(d)$. On same lines of (Delange-style) proof of 
Wintner-Delange formula, we calculate: 
$$
\Wintner_q(F)=\sum_{{{d=1}\atop {d\equiv 0\bmod q}}}^{\infty}{{F'(d)}\over d}
=\sum_{m=1}^{\infty}\sum_{K=1}^{\infty}\mu(K)\widehat{F}(qmK)
=\sum_{K=1}^{\infty}\sum_{{{n=1}\atop {n\equiv 0\bmod K}}}^{\infty}\mu(K)\widehat{F}(qn), 
$$
\par
\noindent
exchanging $m,K$ and, then, applying 1st eq. in $(6)$, a kind of M\"obius inversion, so 
$$
\Wintner_q(F)=\sum_{n=1}^{\infty}\widehat{F}(qn)\sum_{K|n}\mu(K)
=\widehat{F}(q), 
\leqno{(10)}
$$
\par
\noindent
in which, say, we are again exchanging series; this is possible, from absolute convergence of double series in $(10)$, we obtain by $\omega(qn)\ge \omega(n)$ $\Rightarrow $ $2^{\omega(qn)}\ge 2^{\omega(n)}$, uniformly in $q\in \N$: 
$$
\sum_{n=1}^{\infty}\left|\widehat{F}(qn)\right|\sum_{K|n}\mu^2(K)=\sum_{n=1}^{\infty}2^{\omega(n)}\left|\widehat{F}(qn)\right|
\le \sum_{n=1}^{\infty}2^{\omega(qn)}\left|\widehat{F}(qn)\right|
=\sum_{{{m=1}\atop {q|m}}}^{\infty}2^{\omega(m)}\left|\widehat{F}(m)\right|
$$
$$
\le \sum_{m=1}^{\infty}2^{\omega(m)}\left|\widehat{F}(m)\right|, 
$$
\par
\noindent
this last converging absolutely, from our hypothesis $(7)$ above. 
Thus $\widehat{F}=\Wintner(F)$ is proved. 
We still have to prove $\Wintner_q(F)=\Carmichael_q(F)$, $\forall q\in \N$, but this follows from 
$(6)$ and this last bound, in case $q=1$ : they imply the Wintner assumption $(WA)$, 
whence Wintner's Criterion (Theorem 2) proves $\Wintner(F)=\Carmichael(F)$.\hfill $\square$ 

\bigskip

\par
\noindent {\bf 4. Hildebrand and Lucht Lands: Harvesting from Theories.}
\enspace 
We start with a very general result, namely Hildebrand 1984 Theorem [Hi], that ensures existence and pointwise convergence of a R.e. for {\it any} arithmetic function! Even more, the expansion is {\it finite}, actually! 
\par
So, we are quoting {Hildebrand's Theorem}. 
\smallskip 
\par
\noindent
{\bf Theorem 5.}
\enspace 
{\it Let $F:\N \rightarrow \C$. Then, it has a finite Ramanujan expansion }
$$
F(n)=\sum_{q}\widehat{F}(q,n)c_q(n), 
$$
\par
\noindent
{\it $\forall n\in \N$ fixed, with Ramanujan coefficients $\widehat{F}(q,n)\in \C$, eventually depending also on $n$. } 

\medskip

\par
\noindent
Hildebrand proves it exhibiting recursively defined coefficients, a-priori depending on $n$, too [ScSp, p.167]. 
\par
Of course a trivial set of coefficients is, $\forall n\in \N$, 
$$
\widehat{F}(1,n)=F(n),
\quad 
\widehat{F}(q,n)=0,
\quad \forall q>1,
$$
\par
\noindent
but this expansion has only one coefficient ! Another example of finite Ramanujan expansion, depending on $n$, follows writing down $F$ in terms of its Eratosthenes transform: 
$$
F(n)=\sum_{d|n}F'(d)
=\sum_{d\le n}{{F'(d)}\over d}\sum_{\ell|d}c_{\ell}(n)
=\sum_{\ell \le n}\left(\sum_{{d\le n}\atop {d\equiv 0(\!\!\bmod \ell)}}{{F'(d)}\over d}\right)c_{\ell}(n)
=\sum_{\ell \le n}\widehat{F}(\ell,n)c_{\ell}(n),
$$
\par				
\noindent
applying $(2)$ for the condition $d|n$, with a kind of Wintner \lq \lq truncated\rq \rq \thinspace coefficients (where Eratosthenes transform is truncated, there's dependence on $n$): 
$$
\widehat{F}(\ell,n):=\sum_{{d\le n}\atop {d\equiv 0(\!\!\bmod \ell)}}{{F'(d)}\over d}. 
$$ 
\par
\noindent
This may be called the \lq \lq {\it standard finite Ramanujan expansion}\rq \rq, of any $F:\N \rightarrow \C$. 

\medskip

The main problem, in the Theorem, is due to the $n-$dependence of coefficients. The $q-$sum is finite, but again may depend on $n$ ! However, we have {\it possible} {\bf non-uniqueness} of these coefficients (for $F(n)=\0(n)$, $\widehat{F}(q)={1\over q}$ has $(C)$ above, that is clearly not a finite one, but with $\widehat{F}(q)=\0(q)$ it's trivially finite!). 

\medskip

\par
A very easy and non-technical result, we now quote and prove very quickly in three lines, comes from {\it Lucht's Theorem}, Theorem 3.1 [Lu]. 
\smallskip
\par
\noindent
{\bf Theorem 6.}
\enspace 
{\it Let $\widehat{F}:\N \rightarrow \C$ be such that }
$$
d\sum_{K=1}^{\infty}\widehat{F}(dK)\mu(K)
$$
\par
\noindent
{\it converges $\forall d\in \N$. Then }
$$
\sum_{q=1}^{\infty}\widehat{F}(q)c_q(a)=\sum_{d|a}d\sum_{K=1}^{\infty}\widehat{F}(dK)\mu(K) 
$$
\par
\noindent
{\it converges $\forall a\in \N$, to that function of $a$. } 
\par
\noindent {\bf Proof}. Letting $x\to \infty$ may assume $x\ge a$ : apply 2nd equation of Lemma 1, 
$$
\sum_{q\le x}\widehat{F}(q)c_q(a)=\sum_{d|a}d\sum_{{{q\le x}\atop {q\equiv 0\bmod d}}}\widehat{F}(q)\mu(q/d)
=\sum_{d|a}d\sum_{K\le x/d}\widehat{F}(dK)\mu(K), 
$$
\par
\noindent
so the convergence of RHS (we're assuming) implies LHS convergence.\hfill $\square$ 
\medskip
\par
Both Hildebrand \& Lucht Theorems were obtained by means of very powerful theories (see [ScSp], [Hi] and [Lu]), that allowed Lucht [Lu95] to get the Ramanujan coefficients for the $K-$divisor function (compare [CLa] for its numerous links), defined as 
$$
d_K(n)\defineq \multisum_{{{n_1 \enspace \thinspace \cdots \enspace \thinspace \thinspace n_K}\atop {n_1 n_2 \cdots n_K=n}}}1
$$
\par
We quote (a part of) Lucht's Corollary 5.5 in [Lu], recalling: \hfil $p^{\ell}\parallel n\definiz p^{\ell}\vert n, p^{\ell+1}\not\vert n$. 
\smallskip
\par
\noindent
{\bf Theorem 7.}
\enspace 
{\it Fix any integer $K\ge 1$. Then }
$$
\widehat{d}_{K+1}(n)={{(-1)^K}\over {K!}}\;{{\log^K n}\over n} \prod_{p^{\ell}\Vert n}\left(\left(1-{1\over p}\right)^K\sum_{\lambda\ge \ell}{{K+\lambda-1}\choose {K-1}}p^{\ell-\lambda}\right)^{-1}. 
$$

\par
\noindent
(See that case $K=1$ is Ramanujan's [R], while $K\ge 2$ is Lucht's, compare [Lu] \& [Lu95].) 

\bigskip 

\par
\noindent {\bf 5. Finite Ramanujan expansions Ranges: Truncated Divisor Sums.}
\enspace 
Since this is a map, we point to the \lq \lq land\rq \rq, where the {\it finite Ramanujan expansions} (f.R.e.) 
$$
F(n)=\sum_{q\le Q}\widehat{F}(q)c_q(n),
\enspace \forall n\in \N 
\leqno{(FRE)}
$$
\par				
\noindent
are born, i.e., [CMS], [CM2], for the moment (with Ram Murty we expect to supply other papers in the series). So, all things we'll say in this section, actually, (even if we forget) are all quoted from these(by now) two papers. (Recall that, say, Hildebrand's f.R.e. will not be studied in this section.) 
\par
We start with a V.I.P.: 
\smallskip 
\par
\centerline{ALL PURE f.R.e. ARE TRUNCATED DIVISOR SUMS \& VICE VERSA} 
\par
\noindent
namely, next Proposition 2. Recall that we abbreviate t.d.s.$=$truncated divisor sum: 
$$
F:\N \rightarrow \C \; \hbox{\rm is a t.d.s.} 
\; \definiz \; 
\exists F':\N \rightarrow \C, \; \exists Q\in \N \; : \; F(n)=\sum_{d|n,d\le Q}F'(d), \enspace \forall n\in \N
$$ 
\par
\noindent 
where both $Q$ and $F'$ DO NOT depend on $n$, of course. For example \thinspace $d(n)\defineq \sum_{d|n}1$, the divisor function ($=$number of positive divisors of $n$, seen in Theorem 3 proof), is NOT a truncated divisor sum (not only above representation has to be $n-$independent, as stated, but also true \enspace $\forall n\in\N$ : otherwise, $n\le Q$ $\Rightarrow$ $d\le Q$) ! 

\medskip

\par
\noindent
{\bf Remark 2.}
\enspace 
{\it The number $Q\in \N$ is not unique, since we may choose, say, $F'(Q)=0$. Thus, in order to define precisely the {\it range} of a t.d.s., we have to assume $F'(Q)\neq 0$. However, \lq \lq $F$ {\it is of range} $Q$\rq \rq, in our jargon, may also mean: divisors are $d\le Q$, namely $F'(q)=0$, $\forall q>Q$, and $Q$ may not be uniquely defined! } 

\medskip

\par
The following Proposition will keep out of our study, actually, all the f.R.e. that are not pure : sorry for the beautiful Theorem by Hildebrand ! In fact, we know that his f.R.e. are not necessarily pure, because they are valid $\forall F:\N \rightarrow \C$, while the t.d.s. are, say, \lq \lq a drop in this ocean\rq \rq, of all arithmetic functions. 

\medskip

Notice that above we expressed $F(n)$ in TWO WAYS. If $F$ is a t.d.s, through its Eratosthenes Transform $F'$ and through its {\it finite Ramanujan coefficients} $\widehat{F}$, whenever it's a pure f.R.e. This duality $F' \leftrightarrow \widehat{F}$ is perfect, whenever $F$ is AT THE SAME TIME a t.d.s. and a pure f.R.e. : this is proved by next Proposition, that renders crystal clear the duality by formul\ae, for the stated link. 

\bigskip

\par
We are ready to prove, in half page, the following elementary result. 
\smallskip 
\par
\noindent
{\bf Proposition 2.}
{\it Take any $F:\N \rightarrow \C$. Then } 
\par
\centerline{$F$ {\it is a t.d.s. $\Longleftrightarrow$ $F$ has a pure f.R.e.}}
\par
\noindent {\bf Proof}. For $\Rightarrow$, use $(2)$ : 
$$
F(n)=\sum_{{d|n\atop d\le Q}}F'(d)
=\sum_{d\le Q}{{F'(d)}\over d}\sum_{q|d}c_q(n)
=\sum_{q\le Q}\sum_{{{d\le Q}\atop {d\equiv 0\bmod q}}}{{F'(d)}\over d}c_q(n)
=\sum_{q\le Q}\widehat{F}(q)c_q(n), 
$$
\par 
\noindent 
for the {\it finite Ramanujan coefficients} (or $Q-${\it truncated} R.c.s) 
$$
\widehat{F}(q)\defineq \sum_{{{d\le Q}\atop {d\equiv 0\bmod q}}}{{F'(d)}\over d} 
\enspace (\Rightarrow \enspace 
q\le Q, \enspace \hbox{\rm otherwise } \widehat{F}(q)=0). 
\leqno{(Fhat)}
$$
\par
\noindent 
Vice versa, for $\Leftarrow$, use Lemma 1 second equation : 
$$
F(n)=\sum_{q\le Q}\widehat{F}(q)c_q(n)
=\sum_{d|n}d\sum_{{{q\le Q}\atop {q\equiv 0\bmod d}}}\widehat{F}(q)\mu(q/d)
=\sum_{d|n,d\le Q}F'(d), 
$$
\par
\noindent 
for the $Q-$truncated Eratosthenes transform 
$$
F'(d)\defineq d\sum_{K\le Q/d}\widehat{F}(dK)\mu(K) 
\enspace (\Rightarrow \enspace
d\le Q, \enspace \hbox{\rm otherwise } F'(d)=0).
\leqno{(F')}
$$
\par
\noindent
Notice similarity with $(6)$ : actually, this is a finite version of it.\hfill $\square$ 

\medskip
From the begin to the end, the proof assumes the \lq \lq $n-${\bf independence}\rq \rq, for t.d.s. \& for pure f.R.e. (depending on $n$, resp., ONLY in the \lq \lq $d|n$\rq \rq \thinspace and ONLY in the Ramanujan sum $c_q(n)$.) In the following, {\bf this} kind of independence will be implicit for t.d.s. \& pure f.R.e. ! (Compare the concept of \lq \lq fair\rq \rq, in $\S6$.) 

\medskip

\par
Very good news for uniqueness of R.e.s, in following result. 
\smallskip
\par
\noindent
{\bf Corollary 1.}
{\it Completely uniform R.e. have unique coefficients. In particular, pure f.R.e. enjoy uniqueness.} 
\smallskip
Now we will not prove Corollary 1, since it follows both from the Theorem 3 in $\S3$ and from previous Proposition 2. (The interested reader may prove its second part by Carmichael formula, following step-by-step Theorem 3 Proof in the finite part: the $q\le Q$ sum.) 
\par 
Operatively speaking, Corollary 1 tells us to apply the Carmichael Formula to pure f.R.e. (as our most powerful tool). This answers one possible question: why do we care to write t.d.s. as finite Ramanujan expansions ? (We'll also answer \lq \lq on the field\rq \rq, when coming to shift-R.e. in $\S6$.) 
\par
The equation $(FHat)$ in Proposition 2 proof defines $\widehat{F}$ in terms of $F'$ and vice versa, $(F')$ defines $F'$ in terms of $\widehat{F}$ (compare [CM2] formul\ae). Both these formul\ae \thinspace are V.I.P.! 
\smallskip
\par
We give a property of finite Ramanujan coefficients: since $F'$ vanishes after $Q$, in the, say, $Q-${\it truncated counterpart}, $F_Q(n)\defineq \sum_{d|n,d\le Q}F'(d)$, of a given arithmetic function $F=F'\ast \1$, then $\widehat{F_Q}=0$ after $Q$, too and (from above $(FHat)$, compare [CM2], $\S3$) 
$$
{Q\over 2}<q\le Q
\quad \Longrightarrow \quad 
\widehat{F_Q}(q)={{F'(q)}\over q}. 
\leqno{(H)}
$$
\par
\noindent
We quote, from [CM2], this property $(H)$, since it is a V.I.P., for finite Ramanujan expansions: it distinguishes finite R.e. from, say, classical R.e. ! In fact, the operation of truncating the Eratosthenes transform $F'$ (by previous formul\ae, this is equivalent to truncating R.c.s), so to speak, reflects on Ramanujan coefficients in the \lq \lq change\rq \rq, say, of last R.c.s, that we'll call \lq \lq {\bf High} coefficients\rq \rq. (A way to say with high indices, here.) 
\par
While $F$ always has R.c.s (Hildebrand's Theorem) and these, in general, will not vanish definitely (think about $\widehat{\0}_{\rm Ram}(q)=1/q$, esp.), instead, $Q-$truncated R.c.s all vanish after $Q$; so they have in some sense to cope with the original expansion (of $F$, not of its $Q-$truncated counterpart) and this \lq \lq forces\rq \rq, so to speak, the high coefficients to rearrange and recover the final expansion (now, a finite one!), substituting all of the \lq \lq missing R.c.s\rq \rq, with $q>Q$, i.e. the \lq \lq tail\rq \rq. The property above is a V.I.P. just for this reason; but, also, because a simple explicit formula is given for high R.c.s, meaning with $Q/2<q\le Q$ hereafter, in finite expansions. 

\bigskip

\par
Another V.I.P. for R.c.s, which is now a {\bf heuristic property}, is the following: 
$$
q\le Q_0,\enspace \hbox{\rm with } Q_0 \thinspace \hbox{\rm \lq \lq small\rq \rq, w.r.t.}\enspace Q
\quad \Longrightarrow \quad 
\widehat{F_Q}(q)\sim \widehat{F}(q), 
\leqno{(L)}
$$
\par
\noindent
where the $F_Q(n)=\sum_{q|n,q\le Q}F'(q)$ is the $Q-$truncated counterpart of our $F$ (seen above) and, as $Q\to \infty$, the new parameter $Q_0=o(Q)$ is suitably small (esp., think about $Q_0=\sqrt Q$) : say, here these \lq \lq {\bf Low} coefficients\rq \rq, in $(L)$, are asymptotic to classic ones (for $F$). Compare [C] for an application to $F=\Lambda$, von Mangoldt function (for primes). Also, more in general, for truncations of arithmetic functions, see $\S4$ in [CM2]. 
\medskip
\par 
The main question, after knowing the behavior of low \& high finite R.c.s, is of course: what about intermediate ones ? This is a very difficult question to answer ! 
\medskip
\par
Address, once again, for all other f.R.e. features is in [CMS] and [CM2]. 

\bigskip 

\par
\noindent {\bf 6. Shift Ramanujan expansions Landscape: A classic Panorama.}
\enspace 
The \lq \lq realm\rq \rq, say, of shift-Ramanujan expansions started \lq \lq upon a time\rq \rq (few months) ago [CM2]. Thus a genealogy is there, especially (see $(CC)$ next) Theorem 1 and Corollary 1. Here, we are giving a first \lq \lq landscape\rq \rq, for it; in fact, we plan with Ram Murty to explore it much more, in our future papers of the series \lq \lq Finite Ramanujan expansions and shifted convolution sums of arithmetical functions\rq \rq. 
\par				
Given two arbitrary arithmetic functions $f,g:\N \rightarrow \C$, their {\it correlation} (or {\it shifted convolution sum}) is: 
$$
C_{f,g}(N,a)\defineq \sum_{n\le N}f(n)g(n+a), 
$$
\par
\noindent
which, in turn, is itself an arithmetic function of $a\in \N$, called the {\it shift}. (The $N\in \N$ is the {\it length} of $C_{f,g}$.) Hence, if we expand in terms of $a$, we get a kind of new R.e., say, the {\it shift-Ramanujan expansion} (of $C_{f,g}$) 
$$
C_{f,g}(N,a)=\sum_{\ell=1}^{\infty}\widehat{C_{f,g}}(N,\ell)c_{\ell}(a),
\enspace
\forall a\in \N, 
$$
\par
\noindent
that always exists, thanks, again, to Hildebrand's Theorem (see $\S4$) ! 
\par
We may also expand the single $f$ and $g$ inside it : a big surprise came to us [CMS] when we discovered that, however arbitrary are $f:\N \rightarrow \C$ and $g:\N \rightarrow \C$ in their correlation, we always get two single R.e.s which are both {\bf finite} ! This comes immediately from a simple trick, we see now. 
\par
\noindent
This is the {\bf vital remark} [CMS] : 
$$
C_{f,g}(N,a)=\sum_{n\le N}\sum_{d|n}f'(d)\sum_{q|n+a}g'(q)
=\sum_{d}\sum_{q}f'(d)g'(q)\sum_{
{{{n\le N}\atop {n\equiv 0\bmod d}}\atop {n+a\equiv 0\bmod q}}
}
1
$$
$$
=\sum_{d\le N}\sum_{q\le N+a}f'(d)g'(q)\sum_{{n\le N\atop n\equiv 0\bmod d}\atop {n+a\equiv 0\bmod q}}1
=\sum_{n\le N}\sum_{{d|n\atop d\le N}}f'(d)\sum_{{q|n+a\atop q\le N+a}}g'(q), 
$$
\par
\noindent
once written $f(n)=\sum_{d|n}f'(d)$ and $g(m)=\sum_{q|m}g'(q)$ in terms of their Eratosthenes transforms and observed that $d|n$, $n\le N$ $\Rightarrow $ $d\le N$ ! 
\par
As we see, this trick turns $f$ and $g$ into t.d.s., whence (see Proposition 2) we get the {\it finite} R.e.s for our arithmetic functions $f$ \& $g$ (Proposition 2 for formul\ae), inside $C_{f,g}$: 
$$
C_{f,g}(N,a)=\sum_{d\le N}\sum_{q\le N+a}\widehat{f}(d)\widehat{g}(q)\sum_{n\le N}c_d(n)c_q(n+a) 
$$
\par
\noindent
which is, again, valid $\forall a\in \N$; but, doesn't display any $c_{\ell}(a)$, so it is not a shift-R.e. ! 
\par
The point, here, is that in fact these two single R.e.s may help in finding the {\it shift-Ramanujan coefficients} \thinspace $\widehat{C_{f,g}}(N,\ell)$. 
If we may, say, exchange sums applying $(CF)$, then orthogonality \lq \lq reveals\rq \rq \thinspace them! 
\par
\smallskip
\par
After next definition, we need for this, we see how. 
\smallskip
\par
\noindent
We define a correlation \lq \lq {\it fair}\rq \rq, if the shift-dependence (the way it depends on $a$) is only inside the $g$ argument, i.e., $n+a$ :  there's no other dependence on $a$, neither in $f$ and its support, nor in $g$ and its support. (For example, in [CM2] the function $f_H$ given at last in the Appendix has $C_{f_H,f_H}$ which is not fair, since $a\le H$ implies $f_H$ itself depends on $a$, better, its support does!) 
\par
Soon after defining this, in order to use it, we have first to \lq \lq clean up\rq \rq, say, the range of $g$ : it depends on $a$, see the above. For this, we use the simple idea of \lq \lq cutting\rq \rq, with a small remainder (details in $(1)$ of [C]): 
$$
C_{f,g}(N,a)=\sum_{n\le N}\sum_{{d|n\atop d\le N}}f'(d)\sum_{{q|n+a\atop q\le N+a}}g'(q)
=\sum_{n\le N}\sum_{{d|n\atop d\le N}}f'(d)\sum_{{q|n+a\atop q\le N}}g'(q)+O_{\varepsilon}(N^{\varepsilon}(N+a)^{\varepsilon}a)
$$
$$
=C_{f,g_N}(N,a)+O_{\varepsilon}(N^{\varepsilon}(N+a)^{\varepsilon}a), 
$$
\par
\noindent
whenever, of course, $f,g\EssBdd 1$, i.e., they satisfy Ramanujan Conjecture. Here $g_N(m)\defineq \sum_{q|m,q\le N}g'(q)$ is the $N-$truncated counterpart of $g$. In our jargon, a t.d.s. of range $N$. 
\medskip
\par
Now, we start using the concept of fair correlation.
\medskip
\par				
\noindent
The $C_{f,g_N}$ is fair iff in the formula (from $Q=N$ \& $F=g_N$ in $(FRE)$, at $\S5$ begin) 
$$
C_{f,g_N}(N,a)=\sum_{q\le N}\widehat{g_N}(q)\sum_{n\le N}f(n)c_q(n+a) 
\leqno{(11)}
$$
\par
\noindent 
the only, say, place in which there's $a-$dependence is the Ramanujan sum $c_q(n+a)$. This is vital to calculate Carmichael's coefficients (of course, also for their $\exists$), since in $(11)$ we may exchange all the sums, when our (cut-)correlation $C_{f,g_N}(N,a)$ is fair: 
$$
\Carmichael_{\ell}(C_{f,g_N})=\sum_{q\le N}\widehat{g_N}(q)\sum_{n\le N}f(n){1\over {\varphi(\ell)}}\lim_x {1\over x}\sum_{a\le x}c_q(n+a)c_{\ell}(a)={{\widehat{g_N}(\ell)}\over {\varphi(\ell)}}\sum_{n\le N}f(n)c_{\ell}(n),
\leqno{(CC)}
$$
\par 
\noindent 
where we've applied the {\it orthogonality of Ramanujan sums} (proved by Carmichael [Ca] himself, that's why $(CF)$ bears his name), see [M] for a complete proof: 
$$
\lim_x {1\over x}\sum_{a\le x}c_q(n+a)c_{\ell}(a)=\1_{q=\ell}\,c_{\ell}(n),
\quad \forall \ell,n\in \N. 
$$
\par
\noindent
There's a way, now, to get the shift-R.c.s and passes through a remark (many vital remarks give you a vital theory!), say, {\bf Ramanujan inheritance property} : \lq \lq $g$ is a t.d.s\rq \rq \thinspace is inherited by $C_{f,g_N}$ , now becoming a t.d.s. (i.e., shift's divisors $d|a$ are truncated). 
\par
Last but not least, we have to assume Delange Hypothesis $(DH)$ for our $C_{f,g_N}$: 
$$
\sum_{d}{{2^{\omega(d)}|C_{f,g_N}'(N,d)|}\over d}<\infty. 
\leqno{(\ast \ast)}
$$
\par
\centerline{{\bf In [C] the Theorem and its Corollary hold assuming } $(\ast \ast)$.} 
\medskip
\par
\noindent
Of course, Hardy-Littlewood Conjecture [HL] is proved {\it conditionally} in [C], {\it under} $(\ast \ast)$. 
Is there a way to generalize further the results in [C] ? Namely, can we go from Delange to Wintner assumption, maybe making an even weaker assumption ? (In case of shift-R.e.s this is plausible, but for general R.e.s maybe a shortcut to $(DH)$ is impossible! ) 
\par
All of the things we'll say, after this big exposition on s.R.e., were originated (then generalized to R.e.) from our attempts (at present, not useful) to get, from $(CC)$, whence $\Wintner_q(C_{f,g})=0$, $\forall q>N$ (this, under $(WA)$, is proved) that $C_{f,g}'(N,d)=0$, $\forall d>N$ (while, instead, THIS IS HARD) : from which (i.e., $(iv)$, implying $(iii)$ in Theorem 1 of [CM2]) we get next equation (compare $(iii)$ of [CM2] Th.1) called \lq \lq The Reef\rq \rq. 
\par
We wish to derive the \lq \lq {\bf R}amanujan {\bf e}xact {\bf e}xplicit {\bf f}ormula\rq \rq, the {\stampatello Reef} [C], for the (cut-)correlation: 
$$
C_{f,g_N}(N,a)=\sum_{q\le N}\Big({{\widehat{g_N}(q)}\over {\varphi(q)}}\sum_{n\le N}f(n)c_q(n)\Big)c_q(a),
\enspace \forall a\in \N. 
\leqno{{\stampatello Reef}}
$$
\par 
\noindent
This we derived above, using Delange Hypothesis $(\ast \ast)$, in the style of our previous paper [C] and now we wish to obtain some weaker results, under weaker assumptions. (Namely, a weaker version of the Reef, say.) 
\par
See: $(CC)$ works under the only hypothesis of fair correlation (i.e., $C_{f,g_N}(N,a)$ is fair); the problem is that (compare the following), in general, the Carmichael coefficients can't, so to speak, be forced to be Ramanujan coefficients (in other words, the series in $(RE)$ with Carmichael coefficients may, if convergent, not converge to $F(n)$ ! ) and this is our only problem. In fact, thanks to $(CC)$ and the vanishing of $\widehat{g_N}(q)$, after $N$, once given convergence to $C_{f,g_N}(N,a)$ we immediately get the Reef ! 
\medskip
\par
This problem (of convergence of R.e. with Carmichael coefficients) is deep and, so to speak, pervasive in the theory of R.e.! (We may say: ALL of   convergence problems for R.e. in [ScSp] deal with this problem.) Compare [Lu] discussion. 
\par
We'll see, soon after this big exposition for s.R.e., illuminating, general examples. 
\medskip
\par				
\noindent
Thus 
$$
C_{f,g_N}(N,a)=\sum_{{d|a\atop d\le N}}C'_{f,g_N}(N,d)+\sum_{{d|a\atop d>N}}C'_{f,g_N}(N,d), 
\quad
\forall a\in \N,
$$
\par 
\noindent
splitting at $N$ (also, for $g$ of range $Q$, we split at $Q$), after M\"obius inversion [T], from the definition we recall: 
$$
C'_{f,g_N}(N,d)\defineq \sum_{t|d}C_{f,g_N}(N,t)\mu(d/t),
\quad
\forall d\in \N. 
$$
\par
\noindent
From formul\ae, in Proposition 2, we get 
$$
C_{f,g_N}(N,a)=\sum_{q\le N}\widehat{C_{f,g_N}}(N,N,q)c_q(a)+\sum_{{d|a\atop d>N}}C'_{f,g_N}(N,d), 
\quad
\forall a\in \N,
\leqno{(12)}
$$
\par
\noindent
immediately from the definition: 
$$
\widehat{C_{f,g_N}}(N,Q,q)\defineq \sum_{{d\le Q\atop d\equiv 0\bmod q}}{{C'_{f,g_N}(N,d)}\over d}, 
\enspace \forall q\in \N
\leqno{(QRC)}
$$
\par
\noindent
that are, say, a kind of {\it $Q-$truncated} shift-{\it Ramanujan Coefficients}; recall: original ones SHOULD be the same, without truncating $d\le Q$, but this entails convergence problems at once ! (Again, we are \lq \lq playing with Wintner-Delange formula\rq \rq.) Philosophically, $(QRC)$, as $Q\to \infty$, should approximate original shift-Ramanujan coefficients. 

\medskip

\par
The advantage is clear, because now we may apply $(CF)$ to these coefficients (using Corollary 1, since fairness $\Rightarrow$ $\widehat{C_{f,g_N}}(N,Q,q)$ are pure), getting from $(CC)$ and $(12)$
$$
\widehat{C_{f,g_N}}(N,N,q)={{\widehat{g_N}(q)}\over {\varphi(q)}}\sum_{n\le N}f(n)c_q(n)-{1\over {\varphi(q)}}\lim_x {1\over x}\sum_{m\le x}c_q(m)\sum_{{d|m\atop d>N}}C'_{f,g_N}(N,d), 
\enspace \forall q\in \N. 
$$
\par 
\noindent
See that for $q>N$ they vanish, by definition $(QRC)$ above, with $Q=N$. Notice the presence, even in this \lq \lq $d\le Q-$part\rq \rq, of divisors $d>Q$ (set $Q=N$ here) ! Let us abbreviate these limits as 
$$
L(q)=L_{f,g}(q,N)\defineq {1\over {\varphi(q)}}\lim_x {1\over x}\sum_{m\le x}c_q(m)\sum_{{d|m\atop d>N}}C'_{f,g_N}(N,d), 
\enspace \forall q\in \N, 
$$
\par
\noindent
where \thinspace $\exists L(q)\in \C$, $\forall q\in \N$ \thinspace and $\thinspace L(q)=0$ \thinspace after $N$, too. 

\medskip

\par
\noindent
So, only assuming \lq \lq $C_{f,g_N}(N,a)$ is fair\rq \rq, we have proved the {\stampatello Weak Reef} : \enspace $\forall a\in \N$, 
$$
C_{f,g_N}(N,a)=\sum_{q\le N}\left({{\widehat{g_N}(q)}\over {\varphi(q)}}\sum_{n\le N}f(n)c_q(n)-L(q)\right)c_q(a)+\sum_{{d|a\atop d>N}}C'_{f,g_N}(N,d). 
\leqno{{\stampatello (WeakReef)}}
$$

\bigskip

\par
\noindent {\bf 7. Short Neighborhoods: Odds and Ends.}
\enspace 
This Weak Reef, $(WeakReef)$, is a very useful formula, when taking $a-$averages, say \lq \lq short ones\rq \rq, here $a\le A$ with $A\le N$ : as a gift, 
$$
\sum_{a\le A}C_{f,g_N}(N,a)=\sum_{q\le N}\left({{\widehat{g_N}(q)}\over {\varphi(q)}}\sum_{n\le N}f(n)c_q(n)-L(q)\right)\sum_{a\le A}c_q(a), 
\quad \forall A\le N. 
$$
\par 
\noindent
We address the reader, now, to other papers of ours (also future ones), for the s.R.e. and their applications to averages of correlations [CLa]. 
\par				
We come to general R.e., no more ON s.R.e. However, FROM s.R.e., as we say soon before the $(Reef)$. 
\par
\noindent
We saw in $\S6$ the question if it's true that Carmichael coefficients may be taken as Ramanujan coefficients (too general, a complete answer would, even, supersede [ScSp]). 
\par
In order to give a partial answer, we are inspired by [ScSp] in following the \lq \lq {\bf classic mean value}\rq \rq, for a general $F:\N \rightarrow \C$, that might even not exist: it {\bf is simply} $\Carmichael_1(F)$, in our notation. It influences, so to speak, all other Carmichael coefficients (as [ScSp] hints) and, of course, in case $F\ge 0$, this is evident in our next Lemma. 
\smallskip

\par
\noindent
{\bf Lemma 2.}
Let $F(n)\ge 0$, $\forall n\in \N$ and assume $\Carmichael_1(F)=0$. Then $\Carmichael(F)=\0$, i.e., $\Carmichael_q(F)=0$, $\forall q\in \N$. 
\par
\noindent {\bf Proof.} Simply 
$$
{1\over {\varphi(q)}}\Big|{1\over x}\sum_{n\le x}F(n)c_q(n)\Big|\le {1\over x}\sum_{n\le x}F(n), 
$$
\par
\noindent
from the trivial bound on Ramanujan sum and our hypothesis $F\ge 0$; let $x\to \infty$. \hfill $\square$ 

\smallskip

\par
Our Lemma 4 suggests a link between 1st and subsequent Carmichael coefficients: this has been proved by Hildebrand [Hi] already in 1984, for the uniformly-almost-even functions, compare [Lu], Theorem 2.5. 

\smallskip

\par
\noindent
{\bf Remark 3.}
\enspace 
{\it If $F$ satisfies Lemma 4 hypotheses and $F\neq \0$, then, $\Carmichael(F)=\widehat{F}$ implies the absurd $F=\0$. } 

\smallskip
\par
It may seem that $F\ge 0$, not vanishing everywhere and with classic mean-value $0$ is a too strong requirement. Actually, writing the characteristic function of sets $\hbox{\script S}$ as $\1_{\hbox{\script S}}$, whenever the set $\hbox{\script S}\subseteq \N$ has a {\it natural density}, i.e. the limit defined as 
$$
\Carmichael_1(\1_{\hbox{\script S}})=\lim_{x\to \infty}{1\over x}\sum_{n\le x}\1_{\hbox{\script S}}(n)
=\lim_{x\to \infty}{{\#\{n\in \hbox{\script S}:n\le x\}}\over x}
\defineq \delta(\hbox{\script S}) 
$$
\par
\noindent
exists in $[0,1]$, from Lemma 4 (since $\1_{\hbox{\script S}}\in \{0,1\}$ $\Rightarrow$ $\1_{\hbox{\script S}}\ge 0$ and $\1_{\hbox{\script S}}\neq \0$ by $\hbox{\script S}\neq \emptyset$), the 
\smallskip
\par
\centerline{V.I.P.:\hfill non-empty $\hbox{\script S}\subseteq \N$ with $\delta(\hbox{\script S})=0$ have characteristic function with $\Carmichael=\0$.\hfill} 
\smallskip
\par
\noindent
Thus many examples are there ! (Esp., prime numbers have characteristic function with this property, even by \v{C}ebi\v{c}ev bound $\pi(x)\defineq \{p\in \Primes:p\le x\}\ll x/\log x$, [T] : already \thinspace $\pi(x)=o(x)$ suffices for $\Carmichael_1(\1_{\Primes})=0$.) We even have a non-countable family of such sets. 
\smallskip
\par
\noindent
{\bf Remark 4.}
\enspace 
{\it A noteworthy case is the $\hbox{\script S}$ of squares: $\1_{\hbox{\script S}}'=\lambda$, the Liouville function [T] and $\Wintner(\1_{\hbox{\script S}})=\0$ by the Prime Number Theorem [Lu]. In this case, $|\1_{\hbox{\script S}}'|=|\lambda|=\1$, the constant-$1$ function, so the corollary: $\1_{\hbox{\script S}}$ for squares doesn't satisfy $(WA)$. This \thinspace $\1_{\hbox{\script S}}'=\lambda$, \thinspace completely multiplicative, inspires next Lemma 5. } 

\par
We need to define the {\it completely multiplicative} arithmetic functions $f:\N \rightarrow \C$ by $f(ab)=f(a)f(b)$, $\forall a,b\in \N$; in fact, next Lemma uses them: abbreviating \lq \lq c.m.\rq \rq \thinspace for completely multiplicative, from $(4)$, 
$$
F' \enspace \hbox{\rm is c.m.} \enspace \Rightarrow \enspace 
\Wintner_q(F)=\sum_{d\equiv 0\bmod q}{{F'(d)}\over {d}}
={{F'(q)}\over q}\sum_{m=1}^{\infty}{{F'(m)}\over m}
={{F'(q)}\over q}\Wintner_1(F). 
$$
\par
Hence, the following Lemma is proved. 
\smallskip
\par
\noindent
{\bf Lemma 3.}
Let $F:\N \rightarrow \C$ have $F'$ completely multiplicative. Then 
$$
\Wintner_1(F)=0 
\enspace \Rightarrow \enspace 
\Wintner(F)=0 
$$
\par
\noindent
and, assuming $\Wintner_1(F)\neq 0 $ instead, 
$$
\Wintner_q(F)=0, \enspace \forall q>Q
\enspace \Rightarrow \enspace 
F'(q)=0, \enspace \forall q>Q. 
$$

\par
We'll not prove following immediate Lemma. 
\smallskip
\par
\noindent
{\bf Lemma 4.}
Let $F:\N \rightarrow \C$ have $F'\ge 0$. Then 
$$
\Wintner_q(F)=0, \enspace \forall q>Q
\enspace \Rightarrow \enspace 
F'(q)=0, \enspace \forall q>Q. 
$$

\par				
From these two Lemmas, we are led to formulate: 
\smallskip
\par
\noindent {\bf Conjecture 1}. {\it If } $F:\N \rightarrow \C$ {\it has classic mean-value $\Carmichael_1(F)=\Wintner_1(F)\neq 0$, then}
$$
\Wintner_q(F)=0, \enspace \forall q>Q
\enspace \Rightarrow \enspace 
F'(q)=0, \enspace \forall q>Q. 
$$ 
\par
\noindent
Under this Conjecture we may : assume $(WA)$ instead of $(DH)$ in all the results of [C]. 
\par 
In fact, $(CC)$ above gives, assuming $C_{f,g_N}$ is fair, 
$$
\Carmichael_q(C_{f,g_N})=0, \enspace \forall q>N; 
$$
\par
\noindent
then, $(WA)$ for it implies $\Carmichael(C_{f,g_N})=\Wintner(C_{f,g_N})$, so Conjecture 1 for $F(a)=C_{f,g_N}(N,a)$ entails 
$$
\Wintner_q(C_{f,g_N})=0, \enspace \forall q>N
\enspace \Rightarrow \enspace 
C_{f,g_N}'(N,d)=0, \enspace \forall d>N. 
$$
\par
\noindent
This, say finiteness condition, i.e. $(iv)$ in Theorem 1 [CM2], gives the Reef (as $(iii)$ [CM2], Th.1). 
\par
\noindent
For Ramanujan expansions, it seems that a good sign, for R.c.s, is : $\Carmichael(F)=\Wintner(F)$. 
\par
However, we may have $\Carmichael(F)=\Wintner(F)$ outside the Ramanujan Cloud of this $F$ (i.e., they are not $F$ Ramanujan coefficients), as we'll see soon! 
\par 
In fact, a careful analysis of our bound, linking the Carmichael and the Wintner coefficients in above Wintner-Delange proof, drives us towards our following Theorem, we call the Carmichael-Wintner Formula, for this reason; it has the noteworthy consequence of showing $\Carmichael(F)=\Wintner(F)$, when the mean-value of $|F'(d)|$ vanishes : $\Carmichael_1(|F'|)=0$. 
\par
This is our next \lq \lq Slow Decay\rq \rq, say, condition on $F$. In order to prove this Formula, we need, say, a kind of approximate formula, in next Lemma (we prove in half-page). 
\smallskip
\par
\noindent
{\bf Lemma 5.}
{\it Let $F:\N \rightarrow \C$ be any arithmetic function. Then, as $x\to \infty$, for all fixed $q\in \N$, we have} 
$$
{1\over {\varphi(q)}}\cdot {1\over x}\sum_{n\le x}F(n)c_q(n)=\sum_{{d\le x\atop d\equiv 0\bmod q}}{{F'(d)}\over d}+O_q\left({1\over x}\sum_{d\le x}|F'(d)|\right). 
\leqno{(CW)}
$$
\smallskip
The $O-$constant depends ONLY on $q$, here (so $x-$decay isn't affected). 
\smallskip
\par 
\noindent {\bf Proof.} Let's write $F(n)=\sum_{d|n}F'(d)$ as usual, to get 
$$
{1\over {\varphi(q)}}\cdot {1\over x}\sum_{n\le x}F(n)c_q(n)={1\over x}\sum_{d\le x}F'(d){1\over {\varphi(q)}}\sum_{K\le {x\over d}}c_q(dK), 
$$
\par
\noindent
where now in more detail, by exponential sums cancellation seen in Wintner-Delange Formula (Theorem 2) proof,
$$
{1\over {\varphi(q)}}\sum_{K\le {x\over d}}c_q(dK)={1\over {\varphi(q)}}\sum_{j\in \Z^*_q}\sum_{K\le {x\over d}}e_q(jdK)
= 
$$
$$
=\1_{d\equiv 0\bmod q}{x\over d}+O(1)+\sum_{{\ell|d\atop \ell<q}}\1_{(d,q)=\ell}O\left({1\over {\varphi(q)}}\sum_{j\in \Z^*_q}{1\over {\left\Vert{{jd/\ell}\over {q/\ell}}\right\Vert}}\right), 
$$
\par
\noindent
where the $O(1)$ is an absolute constant (coming from : fractional parts are bounded), so this is 
$$
{1\over {\varphi(q)}}\sum_{K\le {x\over d}}c_q(dK)=\1_{d\equiv 0\bmod q}{x\over d}+O(1)+\sum_{{\ell|d\atop \ell<q}}\1_{(d/\ell,q/\ell)=1}O\left({{\ell}\over {\varphi(q)}}\sum_{|j'|\le {q\over {2\ell}}}{{q/\ell}\over {|j'|}}\right)
$$
$$
=\1_{d\equiv 0\bmod q}{x\over d}+O(1)+\sum_{{\ell|d\atop \ell<q}}O\left({q\over {\varphi(q)}}\log (q+1)\right) 
$$
$$				
=\1_{d\equiv 0\bmod q}{x\over d}+O(1)+O\left(q^2 \log (q+1)/\varphi(q)\right)
=\1_{d\equiv 0\bmod q}{x\over d}+O_q(1), 
$$
\par
\noindent
where now the implied constant depends ONLY on $q$, getting from initial formula above 
$$
{1\over {\varphi(q)}}\cdot {1\over x}\sum_{n\le x}F(n)c_q(n)={1\over x}\sum_{d\le x}F'(d)\left(\1_{d\equiv 0\bmod q}{x\over d}+O_q(1)\right), 
$$
\par
\noindent
whence we obtain the thesis.\hfill $\square$ 

This \lq \lq approximate Carmichael-Wintner Formula\rq \rq \thinspace has LHS$\to \Carmichael_q(F)$ and RHS$\to \Wintner_q(F)$ ! 

\par
Now, in order to conclude concordance of Carmichael \& Wintner coefficients in next Theorem, we need to identify the right condition, namely to let the remainder be infinitesimal with $x\to \infty$ : the condition of {\bf Slow Decay}: 
$$
\sum_{d\le x}|F'(d)|=o(x),
\quad \hbox{\rm as} \enspace x\to \infty. 
\leqno{(SD)}
$$
\par
\noindent
This is a way to say: $\Carmichael_1(|F'|)=0$, namely the classic Mean-Value of our $|F'|$ vanishes. 
\smallskip
\par
\noindent
{\bf Remark 5.}
\enspace 
{\it This condition $(SD)$ is a weaker hypothesis, with respect to $(WA)$ : it can be shown (by partial summation, after splitting at $\sqrt x$) that $(WA)$ implies $(SD)$, but (as $F'(d)=1/\log (d+1)$ proves) the converse isn't true! So, our next Theorem is stronger than Wintner's Criterion, i.e., Theorem 2, even if it has a defect, see next Remark. } 
\smallskip
\par
Once given this, we get the \lq \lq Carmichael-Wintner Formula\rq \rq: $\Carmichael(F)=\Wintner(F)$, in our following result (first fix $q\in \N$, then let $x\to \infty$ to get $\Carmichael_q(F)=\Wintner_q(F)$ from $(CW)$ in Lemma 5). Recall Carmichael \& Wintner transforms, from Remark 1. 
\par
The {\bf Carmichael-Wintner Formula} follows. 
\smallskip
\par
\noindent
{\bf Theorem 8.}
\enspace 
{\it Let $F:\N \rightarrow \C$ satisfy $(SD)$. Then $\forall q\in \N$ }
$$
\exists \Carmichael_q(F)\in \C 
\enspace \Longleftrightarrow \enspace 
\exists \Wintner_q(F)\in \C 
\quad \hbox{\it and, if they exist, } \Carmichael_q(F)=\Wintner_q(F). 
$$
\par 
\noindent
{\it Thus $\Carmichael(F)$ exists iff $\Wintner(F)$ exists and, in case they exist, $\Carmichael(F)=\Wintner(F)$.} 

A big defect of our result, w.r.t. the Wintner's Criterion, is that it doesn't get the existence of these coefficients: it has to assume it (maybe, that's why nobody seems to have thought about it before)! However $\exists \Carmichael(C_{f,g_N})$ has, under Conjecture 1, the consequence that $(SD)$ $\Rightarrow $ $(Reef)$, namely Slow Decay $\Rightarrow$ the Reef, compare Remark 5. 

\par
We've found Theorem 8 months ago (generalizing from s.R.e., when preparing the talk given for NTW2017). We didn't yet know of a similar 1987 result of Delange, that we have found only days ago (thanks, Google Scholar) and, while (as you see from Lemma 5 proof) our is elementary, Delange's 1987 Theorem [De87] needs (even if Professor Hubert Delange didn't say explicitly) the zero-free region, for the Riemann zeta-function at least (see His Lemma). In fact, while our previous result can't be applied, for $F'=\lambda$ the Liouville function (see Remark 4: this $F'$ hasn't $(SD)$ above), next Theorem from [De87] can be applied, as $\lambda \ast \1=\1_{\hbox{\script S}}$, the characteristic function of squares (see [T] and Remark 4), is obviously bounded (and so satisfies $(i)$ next): from $\Carmichael(\1_{\hbox{\script S}})=\0$ (see V.I.P. before Remark 4), we get $\Wintner(\1_{\hbox{\script S}})=\0$ (which, of course, is equivalent to PNT and requires a non-trivial arithmetic information: esp., from a zero-free region). 
\par
Next is the Theorem of Delange in [De87]. (However, expressed in our notation.) 
\smallskip
\par
\noindent
{\bf Theorem 9.}
\enspace 
{\it Let $F:\N \rightarrow \C$ be such that, once fixed $q\in \N$, the following two conditions hold: } 
\smallskip
\item{$(i)$} ${\displaystyle \sum_{n\le x}|F(n)|=O(x)}$; 
\smallskip
\item{$(ii)$} $\exists \Carmichael_d(F)\in \C$, $\forall d|q$. 
\smallskip
\par 
\noindent
{\it Then $\exists \Wintner_q(F)\in \C$ and $\Carmichael_q(F)=\Wintner_q(F)$. } 
\smallskip
\par
\noindent
(Actually our $(ii)$ is $(ii)'$, see {\bf 1.5} in [De87], which is equivalent to His $(ii)$, as He says.)
\smallskip
\par
Summarizing the Theorem : if $F$ is bounded on average (condition $(i)$ above), then $\exists \Carmichael(F)$ implies $\exists \Wintner(F)$ and $\Carmichael(F)=\Wintner(F)$ ! So our V.I.P. before Remark 4 also has, consequently, $\Wintner=\0$ ! 
\par				
If we wish to compare condition $(i)$ above with our $(SD)$ condition we might say that they look like \lq \lq independent\rq \rq, meaning that no one seems to imply the other. 
\par
However, considering the proofs, Delange's result (at least, both from this point of view and from previous application to Liouville's function) seems to be much stronger than our. 
\par
Anyway, both these results don't say a word about R.e. of arithmetic functions involved ! It seems, from our proof of Wintner-Delange Formula, that the possibility to exchange sums in the double series in equation $(\ast)$ is the key to prove : Wintner's coefficients of our $F$ are in $\cloudF$ ! So we, as announced soon before the Reef above, still want to \lq \lq jump\rq \rq, from a definitely vanishing $\Wintner(C_{f,g})$, to the Reef; but (even if our attempts and false starts produced,say, all of these odds \& ends) we still haven't found how! A possibility should be to prove Conjecture 1. 

\medskip

\par
Maybe uniqueness of Ramanujan coefficients is a good indication, for our R.e., and our Theorem 3 and Theorem 4 give examples of uniqueness. 

\medskip

\par
\noindent
{\bf Acknowledgments.} Of course, I want to thank again the Organizers of NTW2017 and, also, Ram Murty for our previous common papers, a source of inspiration for [C] and present paper. 

\bigskip

\centerline{\stampatello Bibliography}

\item{[Ca]} R.D. Carmichael, {\sl Expansions of arithmetical functions in infinite series}, Proc. London Math. Society {\bf 34} (1932), 1--26. 

\item{[C]} G. Coppola, {\sl An elementary property of correlations}, arXiv:1709.06445. 

\item{[CLa]} G. Coppola and M. Laporta, {\sl Generations of correlation averages}, J. Numbers Volume 2014 (2014), Article ID 140840, 13 pages http://dx.doi.org/10.1155/2014/140840 (draft, arXiv:1205.1706.) 

\item{[CMS]} G. Coppola, M. Ram Murty and B. Saha, {\sl Finite Ramanujan expansions and shifted convolution sums of arithmetical functions},  J. Number Theory {\bf 174} (2017), 78--92. 

\item{[CM2]} G. Coppola and M. Ram Murty, {\sl Finite Ramanujan expansions and shifted convolution sums of arithmetical functions, II}, J. Number Theory {\bf 185} (2018), 16--47. 

\item{[Da]} H. Davenport, {\sl Multiplicative Number Theory}, 3rd ed., GTM 74, Springer, New York, 2000. 

\item{[De]} H. Delange, {\sl On Ramanujan expansions of certain arithmetical functions}, Acta Arith. {\bf 31}(1976), 259--270. 

\item{[De87]} H. Delange, {\sl On a formula for almost-even arithmetical functions}, Illinois J. Math. {\bf 31} (1987), 24--35. Available online 

\item{[H]} G. H. Hardy, {\sl Note on Ramanujan's trigonometrical function $c_q(n)$ and certain series of arithmetical functions}, 
Proc. Cambridge Phil. Soc., {\bf 20} (1921), 263--271. 

\item{[HL]} G. H. Hardy and J. E. Littlewood, {\sl Some problems of "Partitio numerorum". III: On the expression of
a number as a sum of primes}, Acta Math. {\bf 44} (1923), 1–-70. 

\item{[Hi]} A. Hildebrand, {\sl \"{U}ber die punktweise Konvergenz von Ramanujan-Entwicklungen zahlentheoretischer Funktionen}, Acta Arith. {\bf 44} (1984), 109--140. 

\item{[Lu95]} L. Lucht, {\sl Weighted relationship theorems and Ramanujan expansions}, Acta Arith. {\bf 70} (1995), 25–-42. 

\item{[Lu]} L. Lucht, {\sl A Survey of Ramanujan expansions}, Int. J. Number Theory {\bf 6} (2010), 1785--1799. 

\item{[M]} M. Ram Murty, {\sl Ramanujan series for arithmetical functions}, Hardy-Ramanujan J. {\bf 36} (2013), 21--33. 

\item{[R]} S. Ramanujan, {\sl On certain trigonometrical sums and their application to the theory of numbers}, Transactions Cambr. Phil. Soc. {\bf 22} (1918), 259--276. 

\item{[ScSp]} W. Schwarz and J. Spilker, {\sl Arithmetical functions,  (An introduction to elementary and analytic properties of arithmetic functions and to some of their almost-periodic properties)}, London Mathematical Society Lecture Note Series, {184}, Cambridge University Press, Cambridge, 1994. 

\item{[T]} G. Tenenbaum, {\sl Introduction to Analytic and Probabilistic Number Theory}, Cambridge Studies in Advanced Mathematics, {46}, Cambridge University Press, 1995. 

\item{[W]} A. Wintner, {\sl Eratosthenian averages}, Waverly Press, Baltimore, MD, 1943. 

\bigskip

\leftline{\tt Giovanni Coppola - Universit\`{a} degli Studi di Salerno (affiliation)}
\leftline{\tt Home address : Via Partenio 12 - 83100, Avellino (AV) - ITALY}
\leftline{\tt e-mail : giovanni.coppola@unina.it}
\leftline{\tt e-page : www.giovannicoppola.name}
\leftline{\tt e-site : www.researchgate.net}

\bye